\documentclass[oneside,12pt]{amsart}
\usepackage{amssymb,amscd,amsmath,latexsym}
\usepackage[mathcal]{euscript}
\usepackage{color}

\def\U{{\bf U}}

\begin{document}

\title {A Stringy Product on Twisted Orbifold K-theory}

\author[Alejandro Adem]{Alejandro Adem$^{*}$}
\address{Department of Mathematics
University of British Columbia, Vancouver, B.~C., Canada}
\email{adem@math.ubc.ca}
\thanks{$^{*}$Partially supported by the NSF and NSERC}

\author[Y.~Ruan]{Yongbin Ruan$^{**}$}
\address{Department of Mathematics,
University of Michigan, Ann Arbor, MI 14627}
\email{ruan@math.wisc.edu}
\thanks{$^{**}$Partially supported by the NSF}

\author[B.~Zhang]{Bin Zhang}
\address{Department of Mathematics, Sichuan University}
\email{bzhang@math.sunysb.edu}

\date{\today}

\begin{abstract}
In this paper we define an associative stringy product for the twisted
orbifold $K$--theory of a compact, almost complex orbifold
$\mathcal X$. This product is defined on the twisted $K$--theory
$^\tau K_{orb}(\wedge\mathcal X)$ of the inertia orbifold
$\wedge\mathcal X$, where the twisting gerbe $\tau$ is assumed
to be
in the image of the inverse transgression $H^4(B\mathcal X,\mathbb Z)
\to H^3(B\wedge\mathcal X,\mathbb Z)$.
\end{abstract}

\maketitle
             \def\SS{{\mathbb S}}
         \def\G{{\mathcal G}}
          \def\B{{\mathcal B}}
          \def\H{{\mathcal H}}
          \def\N{{\mathcal N}}
             \def\K{{\mathcal K}}
          \def \x{{\bf x}}
              \def \M{{\mathcal M}}
              \def \C{{\mathbb C}}
          \def\HH{{\mathbb H}}
          \def \Z{{\mathbb Z}}
              \def \R{{\mathbb R}}
              \def \Q{{\mathbb Q}}
              \def \U{{\mathcal U}}
              \def \E{{\mathcal E}}
              \def \z{{\bf z}}
              \def \m{{\bf m}}
             \def \k{{\bf k}}
          \def \n{{\bf n}}
              \def \g{{\bf g}}
              \def \h{{\bf h}}
              \def \V{{\mathcal V}}
              \def \W{{\mathcal W}}
              \def \T{{\mathbb T}}
              \def \X{{\mathcal X}}
          \def \Y{{\mathcal Y}}
          \def \P{{\bf P}}
              \def \F{{\bf F}}
    \def \p{{\mathfrak p}}
    \def \LL{{\mathfrak L}}
    \def \L{{\mathcal L}}
    \def \O{{\mathfrak O}}
    \def \longto{\longrightarrow}
\def \sl{{\frak sl}}
\def \C{{\mathbb C}}
\def \Q{{\mathbb Q}}
\def \vCech{{\v{C}ech\ }}
\def \calP{{\mathcal P}}
\def \H{{\mathcal H}}
\def \Cbar{{\bar{\mathcal C}}}
\def \N{{\mathbb N}}
\def \S{{\Sigma}}
\def \s{{\sigma}}
\def \T{{\mathbb T}}
\def \t{{\tau}}
\def \z{{\bar z}}
\def \bZ{{\bar Z}}
\def \P{{\mathbb P}}
\def \R{{\mathbb R}}
\def \HH{{\mathbb {H}}}
\def \div{{\rm div}}
\def \W{{\mathbb W}}
\def \L{{\mathcal L}}
\def \l{{\lambda}} \def \bl{{\Lambda}}
\def \bu{\bullet}
\def \fd{{\bullet}}
\def \Z{{\mathbb Z}}
\def \e{{\varepsilon}}
\def \ag{{\frak{g}}}
\def \ah{{\frak{h}}}
\def \O{{\mathcal O}}
\def \M{{\mathcal M}}
\def \K{{\mathcal K }}
\def \CC{{\mathcal C}}
\def \XC{{{X}_{\mathcal C}}}
\def \MC{{{\mathcal M}_{\mathcal C}}}
\def \cf {{\mathcal F}}
\def \w{{\wedge}}
\def \x{{\bf x}}
\def \k{{\kappa}}
\def \XCbar{{{X}_{\bar{\mathcal C}}}}
\def \MCbar{{{\mathcal M}_{\bar{\mathcal C}}}}
\def \i{{\sqrt{-1}}}
\def \p{{\partial}}
\def \b{{\delta}}
\def \D{{\Delta}}
\def \G{{\mathcal G}}
\def \o{{\omega}}
\def \g{{\gamma}}
\def \re {\noindent {\it Remark\ \ }}
\def \proof{{\noindent{\it Proof.\ \ }}}

\newcommand{\zb}[1]{\textcolor{blue}{#1}}

\newtheorem{Th}{Theorem}[section]
\newtheorem{cor}[Th]{Corollary}
\newtheorem{lem}[Th]{Lemma}
\newtheorem{prop}[Th]{Proposition}
\newtheorem{claim}[Th]{Claim}

\theoremstyle{definition}
\newtheorem{dfn}[Th]{Definition}
\newtheorem{example}[Th]{Example}
\theoremstyle{remark}
\newtheorem*{remark}{Remark}


\section{Introduction}

    Over the last twenty years, there has been a general trend towards the infusion of physical
    ideas into
    mathematics. One of
    the successful examples in the last few years
    is the subject of twisted K-theory. Interest
in it originates from two different sources
    in physics, the consideration of a D-brane charge on a smooth manifold by
    Witten \cite{W} and the notion of discrete torsion on an orbifold by Vafa \cite{V}.
    In mathematics, there have been important developments connected to this.
    On the one hand, it inspired a
    new subject often referred to as {\em stringy orbifold theory}.
    On the other hand, it revitalized and re-established
    connections to
    many classical topics such as equivariant K-theory, groupoids, stacks and gerbes.
    For smooth manifolds, the mathematical foundation of twisted
    K-theory has been worked out and for any cohomology class
    $\alpha\in H^3(X, \mathbb Z)$,
    one can associate a twisted K-theory $^{\alpha}K(X)$
(see \cite{BM}, \cite{BCMMS}, \cite{K}, \cite{MS}, \cite {AmSg}).
One interesting phenomenon is the difference
    between a torsion class and a non-torsion one: for torsion $\alpha$, we have
    a natural notion of twisted vector bundle or twisted sheaf; for a non-torsion $\alpha$,
    there is no
    geometric notion of vector bundle and one has to use infinite--dimensional analysis.

    The case of an orbifold or even a more general singular space is much
    more interesting. This naturally relates to equivariant theories if we specialize
    to the case of
    $\X=[M/G]$ where $M$ is a smooth manifold and
    $G$ is a compact Lie group acting almost freely on $M$.
    One can consider a cohomology class $\alpha\in H^3(BG,
    \mathbb Z)$ and its corresponding twisted K-theory,
    where $BG$ is the classical classifying
    space for $G$. This was the
    set-up of \cite{AR} for orbifold twisted K-theory
    $^{\alpha}K_{orb}(\X)$ using discrete torsion.
    Twisted K-theory has been generalized to K-theory twisted
    by gerbes  (see \cite{LU}), and also
    using the
    framework of groupoids (see \cite{LTX}).
    In the general case, one can think that the twisting
    is a cohomology
    class\footnote{To be totally precise, we will actually be
    twisting with cocycles.}
    $\alpha\in H^3(B\X, \Z)$ where $B\X$ is now the
    classifying space of the orbifold
    $\X$.

    One advantage of working with
    orbifolds is the nontrivial cohomological counterpart
    called {\em Chen-Ruan} cohomology of orbifolds, $H^*_{CR}(\X,
    \C)$. One can use discrete torsion \cite{R} (more generally
    torsion gerbes \cite{PRY}) to obtain a twisted Chen-Ruan cohomology $H^*_{CR}(\X,
    \L_{\alpha})$. Moreover, the twisted Chen-Ruan cohomology has
    an important internal product making $H^*_{CR}(\X,
    \L_{\alpha})$ a ring. On the other hand, the tensor product produces a
    map
    $$^{\alpha}K_{orb}(\X)\otimes~ ^{\beta}K_{orb}(\X)\rightarrow
    ~^{\alpha+\beta}K_{orb}(\X).$$
    Note that it shifts the twisting to $\alpha+\beta$;
    one natural question
    is if there is an internal "stringy" product for
    $^{\alpha}K_{orb}(\X)$? Freed, Hopkins and
    Teleman
    \cite{FHT}
    have proved the beautiful result that the twisted equivariant
    K-theory $^{\alpha}K^{\rm{dim}~G}_G(G)$ 
for the adjoint action is isomorphic to
    the Verlinde algebra of representations of the central extension of
    the loop algebra $\L G$ for a semi-simple Lie group $G$. This
    algebra carries a very important
    ring structure via the Verlinde product, whose structure constants
    encode
    the information for
    so--called conformal blocks. Using the group structure of $G$,
    one can also construct a ring structure (via the Pontryagin product)
    for $^{\alpha}K^{\rm{dim}~G}_G(G)$; 
    these rings turn out to be
    isomorphic.

    Due to the importance of the Verlinde product in representation
    theory, the existence of a stringy product on the twisted K-theory
    for more general spaces becomes an important question. This is the
    problem we will address in this article and its sequel.

    Our main observation is that there is indeed a stringy product for
    the twisted K-theory of orbifolds.
    Moreover, the
    key information determining such a stringy product does not lie in $H^3(B\X, \Z)$
    as one conventionally believes; instead, it lies in $H^4(B\X,
    \Z)$.
    Given a class $\phi\in
    H^4(B\X, \Z)$, it induces a class $\theta(\phi)\in
    H^3(B\wedge \X, \Z)$ where $\wedge \X$ is the inertia stack of
    $\X$ and thus we can define a twisted K-theory
    $^{\theta(\phi)}K(\wedge \X).$ The inertia stack $\wedge \X$
    can be viewed as the moduli space of constant loops on $\X$. Furthermore, there is a
    key multiplicative formula for $\theta (\phi)$ characterized by
    the effect of $\phi$ on the moduli space $\M$ of
    constant morphisms from a Riemann surface. This map, which can be
    thought of as the inverse of the classical transgression map, appears
    in \cite{DW} for finite group cohomology.
    Based on this we derive a
    simple extension for orbifold groupoids and explicitly prove
    its multiplicative property (a more geometric version of this
    formula appears in \cite{LU2}).

    Our second ingredient is more subtle: experience from Chen-Ruan cohomology
    tells us that
    a naive definition
    does not give an associative product. The reason lies in the fact that the fixed--point
    sets $X_g, X_h$ for
    $g\neq h$ in general do not intersect each other transversely. It is known
    that in Chen-Ruan cohomology
    theory one can correct the naive definition by introducing a certain obstruction
    bundle. Combining these
    two ingredients, we obtain an associative product which can be viewed as
    a K-theoretic counterpart of
    the Chen-Ruan product for orbifold cohomology.

   \begin{Th}
   Let $\mathcal X$ denote a compact, almost complex orbifold, and let
   $\tau$ be a $U(1)$--valued $2$--cocycle for the inertia orbifold
   $\wedge\X$ which is in the image of the inverse transgression.
   Then there is an associative product on $^\tau K_{orb}(\wedge \X)$ which
   generalizes both the Pontryagin and the orbifold cohomology product.
   \end{Th}

    Our construction is in fact motivated by the so--called Pontryagin product
    on $K_G(G)$, for $G$ a finite group, which is what our construction amounts
    to, for $\mathcal X=\wedge [*/G]$ in the untwisted case.
    As an application, we use our construction to clarify the \textsl{twisted Pontryagin
    product}; it may not always exist,
    and when it exists, it may not be unique either.
    We provide an explicit calculation of the inverse transgression map
    for the cohomology of finite groups, showing that in fact it can be computed using
    the natural multiplication
    map $\mathbb Z \times Z_G(h)\to G$, where $Z_G(h)$ denotes the centralizer
    of $h\in G$.
    Using this we exhibit a group, $G=(\mathbb Z/2)^3$ and an integral cohomology
    class $\phi\in H^4(G,\mathbb Z)$ such that under the inverse transgression
    it maps non--trivially
    for every properly twisted sector, yielding an interesting product
    structure on $^{\theta (\phi)} K_G(G)$.

    One of the original motivations for the introduction of the twisted theory
    in orbifolds was the hope
    of describing the cohomology of desingularizations of an orbifold. Joyce
    constructed five classes of
    topologically different desingularization of $T^6/\Z_4$ \cite{J}, arising from
    a representation $\mathbb Z/4 \subset SU(3)$. It is known that Joyce's
    desingularizations are not captured by discrete torsion. For a while, there was the expectation
    that they
    may be captured by 1-gerbes. The computation in \cite{AP}
    shows that the high hopes for
    1-gerbes is probably
    misplaced; however, we notice that $H^4(B(T^6/\Z_4), \Z)$ seems to contain precisely the
    information related to desingularization. We hope to return to this question later.

    We would like to make a comment
    about notation: throughout this paper we will be using the language of
    orbifold groupoids,
    hence given an orbifold $\mathcal X$ we will be thinking of it in terms
    of a Morita equivalence class of orbifold groupoids, represented by $\G$;
    in this context $^\tau K_{orb}(\mathcal X)$ is interpreted as
    $^\tau K(\G)$, using the notion of twisted $K$--theory of groupoids, which
    we will summarize in Section 3.

The results in this article were first announced by the second
author at the Florida Winter School on Mathematics and Physics in
December, 2004. Here, we present our construction for the orbifold
case. The construction for general stacks will appear elsewhere.
During the course of this work, we received an article by
Jarvis-Kaufmann-Kimura \cite {JKK} which also deals with a stringy
product in $K$-theory; indeed the restriction of our twisted
$K$-theory $^\tau K(\wedge \mathcal X)$ to the non-twisted sector
gives their small orbifold $K$-theory $K(|\wedge \mathcal X |)$. The
authors would like to thank MSRI and PIMS for their hospitality
during the preparation of this manuscript, and the third author
would like to thank the MPI--Bonn for its generous support.

    \section{Preliminaries on Orbifolds and Groupoids}

    In this section, we summarize some basic facts about orbifolds, using the
    point of view of groupoids. Our main reference is the
    book \cite{ALR}, but \cite{Moer2} is also a useful introduction.
    Recall that an orbifold structure can be viewed as an orbifold Morita equivalence class of
    orbifold groupoids; we
    shall present all of our constructions in this framework.

    Suppose that $\G=\{s,t: G_1\rightarrow G_0\}$ is an orbifold groupoid, namely, a proper,
\'etale Lie groupoid, we will use $|\G| $ to denote its orbit space,
i.e., the quotient space of $G_0$ under the equivalence relation:
$x\sim y$ iff there is an arrow $g: x\mapsto y$. Conversely, we call
$\G$ an orbifold presentation of $|\G|$.

Recall that a groupoid homomorphism $\phi : \H \to \G$ between (Lie)
groupoids $\H $ and $\G $ consists of two (smooth) maps, $\phi_0:
H_0\to G_0$ and $\phi _1: H_1\to G_1$, that together commute with
all the structure maps for the two groupoids $\H$ and $\G$.
Obviously, a groupoid homomorphism $\phi : \H \to \G$ induces a
continuous map $|\phi |: |\H |\to |\G |$.

\begin {dfn}
Let $\phi$, $\psi: \H \to \G $ be two homomorphisms. A natural
transformation $\alpha$ from $\phi$ to $\psi$ is a smooth map
$\alpha : H_0 \to  G_1$, giving for each $x\in H_0$ an arrow $\alpha
(x): \phi (x) \to  \psi (x)$ in $G_1$, natural in $x$ in the sense
that for any $h : x \to  x'$ in $H_1$, the identity
$\psi(h)\alpha(x)= \alpha (x')\phi (h)$ holds.
\end{dfn}

\begin{dfn}
Let $\phi : \H \to \G$ and $\psi : \K \to \G$ be homomorphisms of
Lie groupoids. The groupoid
fibered product $\H \times _{\G } \K$ is the Lie
groupoid whose objects are triples $(y; g; z)$ where $y\in H_0$,
$z\in K_0$ and $g : \phi (y)\to \psi (z)$ in $G_1$. Arrows $(y; g;
z)\to(y'; g'; z')$ in $\H \times _{\G } \K$ are pairs $(h; k)$ of
arrows, $h : y \to y'$ in $H_1 $ and $k : z \to z'$ in $K_1 $ with
the property that $g'\phi(h) = \psi (k)g$. 
Composition in $\H \times _{\G } \K$ is defined in
the natural way.
\end{dfn}

Next we recall the notion of equivalence of groupoids.

\begin{dfn}
A homomorphism $\phi: \H \to \G$ between Lie groupoids is called an
equivalence if
the map
$t\pi _1: G_1 \ {}_s\times _\phi H_0\to G_0$
is a surjective submersion
and the square
$$
\begin {array}{rcl}
H_1& \stackrel {\phi}{\to}&G_1\\
 (s,t)\downarrow&& \downarrow(s,t) \\
 H_0\times H_0& \stackrel {\phi \times \phi }{\to}&G_0\times G_0
\end{array}
$$
is a fibered product of manifolds.
\end{dfn}

\begin{dfn}
Two orbifold groupoids $\G $ and $\G '$ are said to be orbifold
Morita equivalent if there is a third orbifold groupoid $\H$ and two
equivalences $\phi : \H \to \G$ and $\phi ': \H \to \G '$. An
orbifold homomorphism from $\H$ to $\G$ is a triple $(\K , \epsilon,
\phi)$, where $\K$ is another orbifold groupoid, $\epsilon : \K\to
\H$ is an equivalence and $\phi : \K \to \G $ is a groupoid
homomorphism. The equivalence relation for orbifold homomorphisms is
generated by natural transformations of $\phi$ and equivalences of
$\K$.
\end{dfn}

\begin{dfn}
The category of orbifolds is the category whose objects
are the orbifold Morita equivalence classes of orbifold groupoids and
the
morphisms are equivalence classes of orbifold homomorphisms.
\end{dfn}

\begin{remark}
In this paper, we will use the term homomorphism for a
\textsl{groupoid
homomorphism}, an orbifold homomorphism will be clearly identified
when it arises.
\end{remark}

There are several important constructions which play a fundamental role
in stringy orbifold theory.
 Given $r> 0$ an integer, we can consider the $r$--tuples of
composable arrows in $\G$, i.e.
$$G_r =\{(g_1,\dots, g_r)\in G_1^r ~ | ~ t(g_i)=s(g_{i+1}), ~i=1,\dots , r\}$$
These fit together to form a simplicial space, whose geometric
realization is the
classifying space $B\G$ of the groupoid $\G$. In our discussion of
homological invariants of groupoids, we will be considering cochains
arising from this complex. Recall that the inertia groupoid $\wedge
\G$ is a groupoid canonically associated with $\G$ which is defined
as follows:

    \begin{dfn}
    For any groupoid $\G$, we can associate an inertia groupoid $\wedge \G$  as
    $$(\wedge \G)_0=\{g\in G_1~|~ s(g)=t(g)\}, (\wedge \G)_1=\{(a,v)\in G_2~|~ a\in (\wedge \G )_0\}$$
    where
    $$s(a,v)=a,\,\, t(a,v)=v^{-1}a v.$$
    More generally, we can define the groupoid of $k$-sectors $\G^k$ as
    $$(\G^k)_0=\{(a_1, a_2, \cdots, a_k)\in G_1^k ~|~ s(a_1)=t(a_1)=\cdots=s(a_k)=t(a_k)\}$$
    $$(\G^k)_1=\{(a_1, a_2, \cdots, a_k, u)\in  G^{k+1}_1~|~ s(a_1)=t(a_1)=\cdots=s(a_k)=t(a_k)=s(u)\}$$
    with $s(a_1, \cdots, a_k, u)=(a_1, \cdots, a_k), t(a_1, \cdots, a_k, u)=(u^{-1}a_1 u, \cdots, u^{-1}a_ku)$.
    \end{dfn}
    The construction of the inertia groupoid and $\G^k$ in general is completely functorial. Namely, a homomorphism of
    groupoids
    induces a homomorphism between $k$-sectors and an equivalence of orbifold groupoids induces an
    equivalence betweem them. In case of orbifold groupoids,
    the inertia groupoid can be identified
    as the space of constant loops on $\G$; more generally, $\G^{k-1}$ can be identified as the space of constant
    morphisms from an orbifold sphere with $k$-orbifold points to $\G$. We will come back to these
    descriptions later.

Another important notion is that of quasi--suborbifold; before defining it
we first point out that for a groupoid $\G$,
and an open subset $V\subset G_0$, $\{s,t: s^{-1}(V)\cap t^{-1}(V)\to V\}$ is a groupoid,
which we will denote by $\G |_V$.

    \begin{dfn}
    A homomorphism of orbifold groupoids $\phi: \G\rightarrow \H$
    is a quasi-embedding  if
    \begin{itemize}
    \item
    $\phi: G_0\rightarrow H_0$ is an immersion.
    \item For any $y\in im(\phi)\subset H_0$, with isotropy group $H_y$,
    $\phi ^{-1}(y)$ is in an orbit of $\G $, and for any $x\in \phi
    ^{-1}(y)$, $\phi: G_x \to H_y$ is injective
    \item
    For any $y\in im(\phi)$, and any $x\in \phi
    ^{-1}(y)$,
    there are neighborhoods
    $U_y$ of $y$ and $V_x$ of $x$ such that $\H|_{U_y}=
    U_y\rtimes H_y$,
    $\G |_{V_x}= V_x\rtimes G_x$ and
    $\G |_{\phi ^{-1}(U_y)}\cong (H_y \times _{\phi (G_x)}
    V_x)\rtimes H_y$.
    \item
    $|\phi|: |\G|\rightarrow |\H|$ is
    proper.
    \end{itemize}
    \end{dfn}

    \begin{dfn}
     $\G$ together with $\phi$ is called a \textsl{quasi--suborbifold}  of $\H$.
    \end{dfn}


    The following are important examples of quasi--suborbifolds.
   \begin{example}
    Suppose that $\G=X\rtimes G$ is a global quotient groupoid
    (i.e. a quotient by a finite group). We
    often use the stacky notation $[X/G]$ to denote the groupoid.
An important object is the inertia groupoid $\wedge
\G=(\sqcup_g X_g)\rtimes G$ where $X_g$ is
the fixed point set of $g$ and $G$ acts on $\sqcup_g X_g$ as $h:
X_g\rightarrow X_{hgh^{-1}}$ by $h(x)=hx$. By our definition,
    $\phi: \wedge \G\to \G$
induced by the inclusion map $X_g\rightarrow X$ is a quasi--embedding.
\end{example}

   \begin{example}
    Let $\G$ be the global quotient groupoid defined as in the previous example.
We would like to define an appropriate notion
of the diagonal $\Delta$ for $\G\times \G$. We define it as
    $\Delta=(\sqcup_g \Delta_g)\rtimes (G\times G)$
where $\Delta_g=\{(x,gx), x\in X\}$. Our definition of
quasi--suborbifold includes this example.
\end{example}

    More generally, we define the diagonal $\Delta (\G )$ as the
groupoid fibered product
$\G\times_{\G}\G$.
One can check that $\Delta(\G)=\G\times_{\G}\G$ is locally of the desired
form and hence
a quasi--suborbifold of $\G\times \G$. Notice that the map
$x\to (x, 1_x, x)$ allows us to identify $\G$
as a component of $\Delta (\G)$.

    \begin{example}
    For $l\leq k$, there are natural evaluation morphisms $e_{i_1, \cdots, i_l}: \G^k\rightarrow \G^l$
    given by
    $$e_{i_1, \cdots, i_l}(a_1, \cdots, a_k)=(a_{i_1}, \cdots, a_{i_l}).$$
    Furthermore, we have $e: \G^k\rightarrow \G$ given by
    $$e(a_1, \cdots, a_k)=s(a_1)=t(a_1)=\cdots=s(a_k)=t(a_k).$$
    The latter one corresponds to taking the image of constant morphism.
    We leave as an exercise for the reader to check that
    $e$ and the $e_{i_1, \cdots, i_l}$ are quasi-embeddings and
    that $\G^k$ is a quasi--suborbifold of $\G^l$.
    \end{example}
    One of the main tools is the notion of a normal bundle. If $i: \G\rightarrow \H$ is a
    quasi-embedding, $i^*T\H$ is a groupoid vector bundle over $\G$ such that $T\G$ is a subbundle.
    Then we can define the normal bundle $N_{\G|\H}=i^*T\H/T\G$. $N_{\G|\H}$ behaves as
    the normal
    bundle does for smooth manifolds.

   %
   Next we introduce the notion of intersection for quasi--suborbifolds.

\begin{dfn}
Let $f: \G_1\rightarrow \H, g: \G_2\rightarrow \H$ denote
quasi--suborbifolds. We define their intersection $\G_1\cap\G_2$ as
the restriction of the pullback $\G_1\times_{\H}\G_2$ to the
component $\H$ in $\H\times_{\H}\H$.
\end{dfn}

Note that under this definition it makes sense to intersect a
quasi--suborbifold with itself. Under certain conditions these
intersections can have nice properties, analogous to the situation
for manifolds. We will be interested in the notion of a
\textsl{clean intersection}.

    \begin{dfn}
    Suppose that $f: \G_1\rightarrow \H, g: \G_2\rightarrow \H$ are
    smooth quasi--suborbifolds, we say that $\G_1$ intersects $\G_2$
    \textsl{cleanly}
    if the intersection orbifold $\G_1\cap \G_2$
    is a smooth quasi--suborbifold of $\H$
    (where as before $\H $ is viewed as a component of $\Delta (\H)$)
    such that for every $x\in (\G _1)_0\cap (\G _2)_0$,
    $T_{(x,1_x,x)}(\G_1\cap \G_2)=T_{x} \G_1\cap T_{x} \G_2$.
    \end{dfn}

    \begin{example}
    As we mentioned before, the evaluation map $e: \wedge \G\rightarrow \G$ is
    a quasi--suborbifold. Then, $e$ with itself forms a clean
    intersection. Indeed the question is local, and locally it corresponds
    to the intersection of fixed point sets $V^g\cap V^h$. This is
    clearly a clean intersection.
    More generally,
     $e_{i_1, \cdots,
    i_l}: \G^k\rightarrow \G^l$ is an quasi-embedding. Then,
    two different quasi-embeddings $e_{i_1, \cdots, i_l}, e_{j_1, \cdots,
    j_l}: \G^k\rightarrow \G^l$ form a clean intersection. We
    leave it as an exercise for our readers.
    \end{example}

As in manifold theory, there is also the notion of transversality
for quasi--suborbifolds.

  \begin{dfn}
    Suppose that $f: \G_1\rightarrow \H, g: \G_2\rightarrow \H$ are smooth homomorphisms.
    We say that $f\times g$ is transverse to the diagonal $\Delta$ if locally $f\times g$
    is transverse to every component
    of the diagonal $\Delta$ on the object level. We say that $f, g$ are transverse
    to each other if  $f\times g$ is transverse to the diagonal $\Delta$.
   \end{dfn}

    \begin{example}
    Suppose that $f: \G_1\rightarrow \H, g: \G_2\rightarrow \H$ are
    quasi--embeddings which are transverse to each other.
    Then the intersection $\G_1\cap \G_2$ is a quasi--suborbifold
    of $\H$.
   \end{example}

Note that a clean intersection need not be transverse, this failure
of transversality plays a role in the definition of orbifold
cohomology and K--theory.

\section{Gerbes and Twisted K--Theory}

We now consider the cohomology and $K$--theory
of orbifold groupoids.

\begin{dfn}
Let $\G$ denote a Lie groupoid, then we define the continuous
$U(1)$-valued k--cochains on $\G$ as
    $$C^k(\G, U(1))=\{\phi: G_k\rightarrow U(1)~|~ \phi ~ \textrm{is continuous} \}.$$
The differential on this abelian group (using additive
notation) is
defined via
   \small{ $$\delta\phi(g_1, \cdots, g_{k+1})=\phi(g_2, \cdots
    g_{k+1})+\sum^k_{i=1} (-1)^{i}\phi(g_1, \cdots, g_i
    g_{i+1},\cdots, g_{k+1})+(-1)^{k+1}\phi(g_1,\cdots,g_k).$$}
\end{dfn}

By a result due to
Moerdijk \cite{Moer}, if $\G$ is an \'etale groupoid
then the cohomology of this chain complex is
the \vCech cohomology of $B\G$ with coefficients in the sheaf $C(U(1))$
of $U(1)$-valued
continuous functions over the
classifying space $B\G$.
By the exact sequence
$$ 0\rightarrow \Z\rightarrow C(\R)\rightarrow
C(U(1))\rightarrow 0,$$
we obtain a long exact sequence in cohomology,
$$H^k(B\G, C(\R))\rightarrow H^k(B\G, C(U(1)))\rightarrow
H^{k+1}(B\G, \Z)\rightarrow
    H^{k+1}(B\G, C(\R)).$$
     Since $C(\R)$ is a fine sheaf, the connecting homomorphism is an
isomorphism, and so for $k>0$,
$$H^k(B\G, C(U(1)))\cong H^{k+1}(B\G, \mathbb Z).$$

We recall

    \begin{dfn}
    An $n$-gerbe on $\G$ is a pair $(\H, \theta)$ consisting of
\begin{itemize}
\item
a refinement $\G\stackrel{\epsilon}{\leftarrow} \H$ (i.e., $\epsilon$ is an equivalence)
\item an $(n+1)$--cocycle $\phi: H_{n+1}\rightarrow U(1)$.
\end{itemize}
\end{dfn}

Next we define equivalence of gerbes.

\begin{dfn}
Given two $n$--gerbes $(\H,\theta)$ and $(\H',\theta')$ on $\G$ we have
the following:

\begin{itemize}
\item
$(\H,\theta)$ is equivalent to $(\H',\theta')$ if there is a
common refinement $\H'\leftarrow \H''\rightarrow \H$ such
that induced $(n+1)$--cocycles on $\H''$ are the same.

\item
$(\H, \theta)$ is isomorphic to $(\H', \theta')$ if there is a
common refinement $\H''$ such that the induced
$(n+1)$--cocycles on $\H''$ differ by a coboundary.
\end{itemize}
\end{dfn}

Let $\K\stackrel{\epsilon}{\to} \G$ be an equivalence. If
$(\H,\theta)$ is a $n$-gerbe over $\K$, by definition, it is a
$n$-gerbe over $\G$. For an $n$-gerbe $(\H,\theta)$ over $\G$, $\H
'=\H \times _{\G }\K$ is an orbifold which is equivalent to both
$\H$ and $\K$ by projections $p_1: \H '\to \H$ and $p_2: \H '\to
\K$. So we have the pull-back $n$-gerbe $(\H ', p_1^*\theta)$ over
$\K $. It is easy to see that under these operations, equivalent
(isomorphic) gerbes go to equivalent (isomorphic) ones, therefore
gerbes behave well under orbifold Morita equivalence.

\begin {dfn} An $n$-gerbe on an orbifold is an equivalence class of pairs $(\G , \theta)$, where $\G $ is a
presentation of the orbifold, $\theta $ is a $n+1$-cocycle on $\G $,
and the equivalence relation is gerbe isomorphism.
\end {dfn}

From the definition, it is clear that an $n$-gerbe defines a \vCech
$(n+1)$-cocycle for the sheaf $C(U(1))$ of continuous $U(1)$--valued
functions on the classifying space $B\G$ ($B\H $ and $B\G $ are
weakly homotopy equivalent). Hence it will define a cohomology class
in $H^{n+1}(B\G, C(U(1)))\cong H^{n+2}(B\G, \Z).$ The image of
$\theta$ under the connecting homomorphism in $H^{n+2}(B\G, \Z)$ is
called its characteristic class or Dixmier-Douady class.


We can associate twisted $K$-theory to a 1-gerbe. For simplicity
    we assume that the $2$--cocycle $\theta$ is defined on $\G$, i.e., we are dealing with
    $(\G, \theta )$, and the twisted $K$-theory $^\theta K (\G )$ will be defined. We follow the
treatment of \cite{LTX}
    to describe this.
    Let $\HH$ be a separable Hilbert space; it is well-known that
    the characteristic class of a principal $PU(\HH)$--bundle
    over $\G$ also lies in $H^3(B\G, \Z)$.
    Hence, given a 1-gerbe,
    we should be able to associate a $PU(\HH)$ bundle with the same characteristic class;
    in fact we can associate
    a canonical principal $PU(\HH)$--bundle.
    We outline its construction.

    For the orbifold groupoid $\G=\{s,t: G_1\rightarrow G_0\}$, let $R=G_1\times U(1)$ be the topologically trivial central extension, and
     $$(g_1,r_1)(g_2, r_2)=(g_1g_2, \theta(g_1, g_2)r_1r_2), $$
     which makes $\{\tilde{s},\tilde {t}: R\rightarrow G_0\}$
     a Lie groupoid, where $\tilde {s}(g,r)=s(g), \tilde {t}(g,r)=t(g)$.

     Now let $G ^x=t^{-1}(x)$; there is a system
     of measure (Haar system) $\lambda=(\lambda^x)_{x\in G_0}$, where $\lambda^x$ is  a measure
     with support $G^x$ such that for all $f\in C_c(G_1)$, $x\rightarrow \int_{g\in R^x} f(g)\lambda^x(dg)$ is
     continuous.
     By $\L^2_x$, we denote the space $L^2(G^x)$ consisting of functions defined on $G^x$
     which are $L^2$ with respect to the Haar measure. Let
     $E_x=\L^2_x\otimes {\HH}$, $E= \sqcup_x E_x.$
     Then $E$ is a countably generated continuous field of infinite dimensional Hilbert spaces over
     the finite dimensional space $G_0$, and therefore is a locally trivial Hilbert bundle according to
     the Dixmier-Douady theorem \cite {DD}.

     The Lie groupoid $R$ acts naturally on $E$: $U(1)$ acts on $\HH$ by complex multiplication.
     Therefore, $E$ is naturally a Hilbert bundle over $\{\tilde {s},\tilde {t}: R\rightarrow G_0\}$.
     Notice $E$ is not a Hilbert bundle over $\G$.
     However, $P(E)$ is a projective bundle over $\G$ with precisely the same characteristic class of $\theta$.
     Let $\B$ be the principal bundle of orthonormal frames of $E$; it is a $U(\HH)$-principal bundle.
     By our previous argument, $P\B$ is a principal $PU(\HH)$-bundle over $\G$.

     Let $\theta$ be a 1-gerbe and $P_{\theta}$ be the associated $PU(\HH)$-bundle constructed above.
     Let $Fred^0(\HH)$ be the space
     of Fredholm operators endowed with the $*$-strong topology and $Fred^1(\HH)$ be the space of self-adjoint
     elements in $Fred^0(\HH)$. Let $\K(\HH)$ be the space of compact
     operators
     endowed with the norm-topology. Now consider the associated bundles
     $$Fred^i_{\theta}(\HH):=P_{\theta}\times_{PU(\HH)} Fred^i(\HH)\rightarrow G_0,$$
     $$\K_{\theta}(\HH):=P_{\theta}\times_{PU(\HH)} \K(\HH)\rightarrow G_0.$$
    By $\F^i_{\theta}$, we denote the space of norm-bounded, $G_1$-invariant, continuous sections $x\rightarrow T_x$
    of the bundle $Fred^i_{\theta}(\HH)\rightarrow G_0$ such that there exists a norm-bounded, $G_1$-invariant,
    continuous section $x\rightarrow S_x$ of $\K_{\theta}(\HH )\rightarrow G_0$ with the property that $1-T_xS_x$ and
    $1-S_xT_x$ are continuous sections of $\K_{\theta}(\HH)$ vanishing at infinity of $|\G |$.
    \begin{dfn}
     For any section $T$ of $\F^i_{\theta}$. We define the support $supp(T)$ as the set of point
     $x\in |\G|$ such that $T_{x'}$ is not invertible for any $x'$ in the preimage of $x$.
    \end{dfn}
    Then we have
    \begin{dfn}
    Let $\G$ be an orbifold groupoid and $(\G ,\theta)$ be a 1-gerbe. We define its
    $\theta$--twisted K-theory as
    $$^{\theta}K^i(\G)=\{[T]~|~T\in \F^i_{\theta}\},$$
    where $[T]$ denotes the homotopy class of $T$ where $T$ is compactly supported.
    \end{dfn}

    Note that since the space of invertible operators is contractible,
    any $T$ with compact support is homotopic to
    a section which is the identity outside a compact subset. Suppose that $i: U\to G_0$ is an open subset,
    using the property above, we have a natural extension
    $$i_*:~^{i^*\theta} K^i(\G|_U)\rightarrow ^{\theta}K^i(\G).$$

\begin{remark}
Suppose that we have a cocycle $\alpha=\beta+\delta \rho$.
    Then there is a canonical isomorphism between central extensions of
    groupoids
    $$\psi_{\rho}: R_{\alpha}\rightarrow R_{\beta}$$
    given by by $\psi_{\rho}(g,r)=(g,\rho(g)r).$
    Hence it induces an isomorphism
    $$P_{\alpha}(E)\rightarrow P_{\beta}(E)$$
    and also a canonical isomorphism
    $$\psi_{\rho}:~^{\alpha}K^i(\G)\rightarrow~ ^{\beta}K^i(\G).$$
    Suppose that in fact $\rho$ is a cocycle, i.e., $\delta \rho=0$.
    Then $\alpha=\alpha+\delta\rho$ and hence we have an
    automorphism $\psi_{\rho}:~^{\alpha}K^i(\G)\rightarrow ^{\alpha}K^i(\G)$.
    Furthermore, if $\rho=\delta \gamma$ is a coboundary,
    then $\psi_{\rho}$ is the identity. Hence $H^1(B\G, U(1))$ acts as
    automorphisms of twisted K-theory.
    It is easy to check in many examples that they are nontrivial
    automorphisms. In the
    literature, twisted K-theory is often referred to as being twisted
    by a \vCech cohomology class or characteristic
    class  of a 1-gerbe. This is a rather ambiguous statement,
    as cohomologous 1-gerbes induce isomorphic twisted
    K-theory, but this is not canonical.
    This observation is particularly important when we define
    a product structure on twisted K-theory.
\end{remark}

To summarize: for 
an orbifold groupoid $\G$, and a 1-gerbe $(\G , \theta)$, twisted $K$-theory
$^\theta K(\G)$ is well--defined up to isomorphism (see \cite{LTX}).

\begin {dfn} For an orbifold and a 1-gerbe on it, taking a presentation
$(\G , \theta)$ of the gerbe, the twisted $K$-theory of the orbifold
is defined to be $^\theta K(\G)$.
\end{dfn}

There is a natural addition operator for 
$^{\alpha}K^*(\G)$ induced by the Hilbert space addition
$\HH\cong\HH\oplus \HH$. On the other hand, the multiplication 
operation induced by Hilbert space
    tensor product $\HH\cong \HH\otimes \HH$ shifts the twisting
$$^{\alpha}K^*(\G)\otimes~^{\beta}K^*(\G)
\rightarrow~^{\alpha + \beta}K^*(\G).$$

\begin {remark} Strictly speaking, there is an issue about canonicity
because of the way we identify $\HH \oplus\HH$ and $\HH\otimes \HH$
with $\HH$ . Since $U(\HH)$ is contractible, any identification will
give the same homotopy classes, therefore the same $K$-theory
element.
\end{remark}
We also want to point out that the product of an element in twisted
$K$-theory with a vector bundle is to be understood in the following
sense. Let $P$ be a family of Fredholm operators on $\HH $
parameterized by a space $M$ and $E$ a complex vector bundle over
$M$ of finite rank. Then $P\cdot E$ is a family of Fredholm
operators on $\HH \otimes E$ parameterized by $M$, which at every
point $x\in M$ has $P\cdot E (x)=P(x)\otimes Id_E$. Hence it is easy
to see that if $P\in ^\theta K^i(\G)$ and $E$ is a $\G$-bundle, then
$P\cdot E \in ^\theta K^i(\G)$ (just like the above remark, the way
to identify $\HH \otimes E$ with $\HH$ does not change the element).
So even though $E$ may not be an element of $K^0(\G)$, the product
with $E$ makes sense.

%

To define our stringy product, we will need a version of the
push-forward map in the context of twisted K-theory. For smooth
manifolds, such a push-forward map has already been worked out by
Carey-Wang \cite{CW}. In the case of orbifold groupoids, we follow
their treatment, the extra effort will be needed to deal with the
groupoid structure.

Let $\G$, $\H $ be almost complex groupoids, $f:\G \to \H$ be a
homomorphism which preserves the almost complex structures, and
$\alpha $ be a $1$-gerbe on $\H$. We will define the push-forward
map $f_*:~^{f^*\alpha} K^*(\G )\to ^{\alpha }K^*(\H )$.

Let $\G=\{s,t: G_1\rightarrow G_0\}$ be a groupoid and $E$ a rank
$n$ complex vector bundle over $\G$, i.e. $\pi : E\to G_0$ is a
complex vector bundle with compatible $\G $-action. We first
establish the Thom homomorphism $\Phi: ~^\alpha K^*(\G )\to ~^{\pi
^*\alpha }K^*(\G \rtimes E)$, where $\G \rtimes E$ is the
transformation groupoid with object set $E$ and arrow set $G_1\times
_{G_0} E$.

Fix any invariant hermitian metric on $E$ \cite {Rj}, then for any
$g \in G_0$, $e\in E_g=\pi ^{-1}(g)$, we use $e^*$ to denote the
dual of $e$ with respect to the fixed hermitian metric.

The complex $\G$-bundle $\pi: E\to G_0$ defines a complex of $\G
$-bundles over $E$,
$$\lambda _E=(\Lambda ^{even} \pi^ *E, \Lambda ^{odd}\pi ^* E , \phi),$$
where $\phi _{(g,e)}=e\wedge -e^*\llcorner$.
Notice that in ordinary $K$-theory, this is the Thom element.

For any element $x\in~ ^\alpha K^*(\G)$, it is represented by a $\G $-equivariant section $x: G_0 \to
Fred^*_\alpha(\HH)$ with $supp(x)$ compact.

\begin {prop} For any $1$-gerbe $\alpha$, there is a $K^0(\G)$-module homomorphism
$$\Phi :~^\alpha K^*(\G )\to ^{\pi^*\alpha}K^*(\G \rtimes E)$$
which is the standard Thom isomorphism in the case of equivariant $K$-theory.
\end{prop}

\proof For any element $x\in~^\alpha K^i(\G)$, it is given by a $\G $-equivariant section $x: G_0\to Fred ^i_\alpha(\HH )$.
Therefore, we have a section $\pi ^*(x): E\to \pi ^* Fred^i_\alpha (\HH)$. Now $supp (\pi
^*(x))=\pi^{-1}(supp(x))$.

For any point $(g, e)\in E$, let us consider the operator:
$$D: \HH \otimes \Lambda ^{even} E_g\oplus \HH \otimes \Lambda ^{odd}E_g \to
\HH \otimes \Lambda ^{even} E_g\oplus \HH \otimes \Lambda ^{odd}E_g$$
defined by
$$D_{(g,e)}=\left (\begin{array}{cc}
(\pi^*x)(g)\otimes 1& -1\otimes \phi_{(g,e)} ^*\\
1\otimes \phi_{(g,e)}&(\pi^*x)(g)^*\otimes 1
\end{array}
 \right)
$$
where $*$ means the adjoint operator. This is the so-called ``graded tensor product" of Fredholm operators.
It is easy to check that $D_{(g,e)}$ is Fredholm and if $e\not= 0$,
then it is invertible.

Globalizing this construction, we have a family $D$ of Fredholm
operators parameterized by $E$. In fact, it is a fiberwise Fredholm
operator on Hilbert bundles $\Lambda ^\bullet \pi^* E\otimes \HH$
over $E$, as we remarked, the identification of fiber with $\HH $
does not matter.

Notice that
\begin {itemize}
\item $supp(D)=\pi^{-1}supp(x)\cap i(G_0)$, where $i:G_0\to
E$ is the zero section.
\item $D$ is a section of $\pi^*Fred^i_\alpha (\HH )\cong
Fred^i_{\pi^*\alpha}(\HH )$.
\item $D$ is $\G \rtimes E$-equivariant.
\end{itemize}
Therefore $D$ represents an element of $^{\pi^*\alpha}K^i(\G \rtimes
E)$

Now we define
$$\Phi :~^\alpha K^*(\G )\to ~^{\pi^*\alpha}K^*(\G \rtimes E)$$
$$x\mapsto D$$
Up to homotopy, as in ordinary $K$-theory, this definition does not
depend on any choice, so it is well-defined. Because it is the
graded tensor product with $\lambda _E$, $\Phi $ is a $K(\G
)$-module homomorphism. Furthermore, by the definition, we see that
it is a generalization of the Thom isomorphism in equivariant
$K$-theory.

We need a slightly more general version of this. Let $U$ be an open
neighborhood of the zero section, from the definition, $\Phi (x)$ is
supported on the zero section, so by restriction, we have the
following Thom homomorphism.

\begin {prop}
For any $1$-gerbe $\alpha$, there is a $K(\G)$-module homomorphism
$$\Phi : ^\alpha K^*(\G )\to ^{\pi^*\alpha}K^*(\G \rtimes E| _U)$$
\end{prop}

To handle the general situation, let us first recall a lemma for Lie
groupoids \cite {Moer0},

\begin {lem} Let $p: F\to G_0$ be a smooth surjective submersion, then the groupoid $F\times _p \G $ is equivalent to
$\G$, where $(F\times _p \G )_1=(F \times F) \ _{p\times p} \times _{s\times t} G_1$ and $ (F\times _p \G )_0=F$, and the
new source map is $s: (F \times F) \ _{p\times p} \times _{s\times t} G_1 \to F, ((x,y),g)\mapsto x$, the new target map
is $t: (F \times F) \ _{p\times p} \times _{s\times t} G_1 \to F, ((x,y),g)\mapsto y$.
\end{lem}

For a Lie groupoid homomorphism $f: \G \to \H $, if we apply the
lemma above to the space $F=G_0\times H_0$, and take $\K$ to be
$F\times _p \H $, we can prove the next lemma, where all the maps
are the natural ones.

\begin {lem} Let $f: \G \to \H $ be a homomorphism of Lie groupoids, then there exists a Lie groupoid
$\K $ and homomorphisms $g: \G \to \K$, $h: \K \to \H$, where
$g_0:G_0\to K_0$ is an embedding and $h$ is an equivalence, such
that $f$ is the composition of $g$ and $h$. In other words, any
homomorphism is an embedding on the object level up to Morita
equivalence.
\end{lem}

Now we can prove the existence of a push-forward map in twisted
$K$-theory.

\begin {Th} If $f: \G \to \H $ is a homomorphism between almost complex groupoids which preserves the almost
complex structures, such that $|f|: |\G|\to |\H|$ is proper, 
and $f_1(\G .x)=\H . f_0(x)$ for any $x\in G_0$,
then there is a push-forward map
$$f_*:~^{f^*\alpha}K^*(\G )\to ~^{\alpha} K^*(\H)$$
\end{Th}

\proof 
Given our last lemma, we may assume that $f_0: G_0\to H_0$ is 
a proper embedding.
By our assumption, the normal bundle $E$ of $G_0$ in $H_0$ is a
complex $\G$-bundle, and we can identify an open neighborhood $U$ of
the zero section in the normal bundle with a neighborhood of
$f_0(G_0)$ in $H_0$, i.e. we have an embedding $j: U\to H_0$ as an
open subset. It defines a homomorphism: $j_*:~^{j^*\alpha} K^*(\G
\rtimes E| _U) \to ~^{\alpha }K^*(\H | _{j(U)})$, because the action
of $\G $ on the normal bundle is induced from the $\H$ action, in
this case any $\G$-equivariant section is $\H $-equivariant.

It is clear that $f\pi$ is homotopic to $j$.
Therefore, $\pi^*f^*\alpha=j^*\alpha+\delta \rho$ for some $\rho$.
The choice of $\rho$ is not unique; for
example, we can add a 1-cocycle. This corresponds exactly to the non-canonicity
of the dependence of twisted K-theory
on the cohomology class of a 1-gerbe. However, $f\pi=j$ on the zero section; therefore,
we can choose $\rho$ such
that $\rho=0$ on the zero section. Since $U$ deformation retracts to the
zero section, it fixes $\rho$ uniquely.

Now we have homomorphisms:
$$^{f^*\alpha} K^*(\G )\stackrel{\Phi }{\longrightarrow}
\ \ ^{\pi ^*f^* \alpha }K^*(\G \rtimes E|_U) \stackrel
{p_{\rho}}{\longrightarrow} ~^{\alpha} K^*(\G\rtimes E|_U) \stackrel
{j_*}{\longrightarrow} \ \ ^{\alpha }K^*(\H |_{j(U)}) \to ~^{\alpha
}K^*(\H) $$ where the last homomorphism is extension for an open
saturated subgroupoid. The composition is the push-forward map
$f_*$.

Given our explicit definition, it is easy to check the following properties of the push-forward map.

\begin {prop} Let $f:\G\to \H$ as before,
then there exists an element $c=c(\G, \H)$ such that for any
$a\in ^{f^*\alpha}K^*(\G)$ and $b\in ^{\beta}K^*(\H )$, we have
$$f^*f_*(a)=a\cdot c$$
$$f_*(a\cdot f^*(b))=f_*(a)\cdot b$$
\end{prop}

\noindent In particular for quasi--suborbifolds, we have following result.

\begin {cor} If $i: \G _1\to \G $ is a quasi--suborbifold, then there is a push-forward map
$$i_*:~^{i^*\alpha}K^*(\G _1)\to ~^\alpha K^*(\G)$$
satisfying the above properties.
\end{cor}

For later purposes we would like to introduce

\begin{dfn}
If $E\to\G$ is a complex orbifold bundle, then its $K$--theoretic
Euler class $e_K(E)$ is defined as $i^*\lambda_E$, the complex of $\G$--vector
bundles obtained by pulling back the Thom element $\lambda_E$ using
the $zero$--section $i:\G\to E$.
\end{dfn}

\noindent Note that we can define the product $x\cdot e_K(E)$ as

$$x\cdot e_K(E) = x\cdot (\wedge^{even}~E) + (x\cdot \wedge^{odd}~E)^*.$$


\section{The Inverse Transgression for Groupoids}

In order to define the stringy product in twisted K--theory, we will
need a cohomological formula to match up the levels which appear in
the twistings. The basic construction is the \textsl{inverse
transgression}, which was defined in \cite{DW}. We provide a
formulation for groupoids inspired by the case of finite groups. We
will also provide some explicit calculations. See \cite{LU} for a
more geometric view on this.

     Recall that
     $$(\wedge \G)_0=\{a\in G_1~|~ s(a)=t(a)\},\,\,
(\wedge\G)_1=\{(a,u_1)\in G_1\times G_1~|~ s(a)=t(a)=s(u_1)\}.$$
     It is easy to check that the $k$--tuples of
composable arrows in $\wedge \G$ are
     $$(\wedge \G)_k=\{(a,u_1,\cdots,u_k)\in
     G_{k+1}~|~ s(a)=t(a)=s(u_1), t(u_i)=s(u_{i+1})\}.$$

     \begin{dfn}
     Define $\theta: C^{k+1}(\G, U(1))\rightarrow C^{k}(\wedge \G,
     U(1))$ by
     $$\theta(\phi)(a, u_1, \cdots, u_k)=(-1)^k\phi(a, u_1, \cdots,
     u_k)+\sum_{i=1}^k (-1)^{i+k}     \phi(u_1, \cdots, u_{i}, a_i, u_{i+1}, \cdots, u_k),$$
     where $a_i=(u_1\cdots u_{i})^{-1} a u_1\cdots u_i.$
     \end{dfn}

A routine but slightly tedious computation shows that this is in fact a
cochain map, i.e. $\delta \theta= \theta\delta.$

We should note that $\theta$ is a natural map defined for all
groupoids. For orbifold groupoids it induces a homomorphism
$$\theta_*: H^k(B\G, U(1))\to H^{k-1}(B\wedge\G, U(1)),$$
and hence a homomorphism
$$\theta_*: H^{k+1}(B\G, \Z)\to H^k(B\wedge\G,
    \Z).$$
The cochain map $\theta$ and the induced map in cohomology
will be called the inverse transgression.


Recall that the moduli space of constant morphisms
$\overline{\M}_3(\G)$ from an orbifold sphere with
three orbifold points can be identified with the $2$-sector orbifold $\G^2$, where
 $$(\G^2)_0=\{(a,b)\in G_2~|~ s(a)=t(a)=s(b)=t(b)\},$$
 $$(\G^2)_k=\{(a, b, u_1,\cdots,u_k)\in G_{k+2}~|~ s(a)=t(a)=s(b)=t(b)=s(u_1), t(u_i)=s(u_{i+1})\}$$
 with
 $$s(a,b,u_1, \cdots, u_k)=(a,b), t(a,b,u_1, \cdots, u_k)=(a_k, b_k)$$
 where
 $$a_i=(u_1 \cdots u_i)^{-1}a u_1\cdots u_i, ~~b_i=(u_1 \cdots u_i)^{-1}b u_1\cdots
    u_i.$$
      There are three natural evaluation morphisms
      $$e_1: \G^2\rightarrow \wedge \G \mbox{ by } e_1(a,b)=a,$$
      $$e_2: \G^2\rightarrow \wedge \G \mbox{ by } e_2(a,b)=b,$$
    $$e_{12}: \G^2\rightarrow \wedge \G \mbox{ by } e_{12}(a,b))=ab.$$
    Furthermore, $e_1,e_2, e_{12}$ are all quasi-embeddings.

    \begin{dfn}
    Define $\mu: C^{k+2}(\G, U(1))\to C^k(\G^2, U(1))$ by
    $$\mu(\phi)(a,b,u_1, \cdots, u_k)=\hskip 2in$$
    \small{$$\phi(a,b,u_1, \cdots,
    u_k)+\sum_{\{(i,j) ~|~ 0\leq i\leq j\leq k ~~(i,j)\ne (0,0)\}}(-1)^{i+j}\phi(u_1, \cdots,
    u_i,a_i,u_{i+1}, \cdots,u_j, b_j,u_{j+1}, \cdots, u_k).$$}
    \end{dfn}

    A second key \textsl{multiplicative formula} is given by the equation

    $$\mu\delta +\delta \mu  ~=~   e^*_1\theta + e^*_2\theta
    -e^*_{12}\theta.$$
    Note that the function $\mu$ defines a chain homotopy between
    $e^*_1\theta + e^*_2\theta$ and $e^*_{12}\theta$.
    If $\phi$ is a cocycle, then
    $\theta(\phi)$ is a cocycle and the formula above implies
    that
    $$e^*_1\theta(\phi)+e^*_2\theta(\phi)=e^*_{12}\theta(\phi)
                      + \delta\mu (\phi).$$
    In particular we see that the difference between
    the cocycles is given by a \textsl{canonical} coboundary,
    expressed explicitly as a function of $\phi$. This will be
    very important when we make our identifications in twisted
    K--theory.

    We will verify and apply this formula in low degree, which is
    our main interest here.


\begin{prop}

Let $\phi$ be an element in $C^3(\G, U(1))$. Then
$$\delta\mu (\phi) + \mu\delta (\phi) = e_1^*\theta (\phi)
+ e_2^*\theta (\phi) - e_{12}^*\theta (\phi).$$

\end{prop}

\begin{proof}
This can be proved by an explicit calculation.

\small{$$\delta\mu (\phi)(a,b,u_1,u_2) = -\phi(a_1,b_1,u_2)+\phi(a,b,u_1u_2)-\phi(a,b,u_1)
+\phi(a_1,u_2,b_2)-\phi(a,u_1u_2,b_2)$$
$$\hskip .5in +\phi(a,u_1,b_1)-\phi(u_2,a_2,b_2)
+\phi(u_1u_2,a_2,b_2)-\phi(u_1,a_1,b_1)$$
$$\mu\delta(\phi)(a,b,u_1,u_2)=~~~\phi(b,u_1,u_2)-\phi(u_1,b_1,u_2)
+\phi(u_1,u_2,b_2)+\phi(a_1,b_1,u_2)
-\phi(a_1,u_2,b_2)$$
$$~~~+\phi(u_2,a_2,b_2)
-\phi(ab,u_1,u_2)+\phi(au_1,b_1,u_2)-\phi(au_1,u_2,b_2)
-\phi(u_1a_1,b_1,u_2)$$
$$~~~+\phi(u_1a_1,u_2,b_2)-\phi(u_1u_2,a_2,b_2)+
\phi(a,bu_1,u_2)-\phi(a,u_1b_1,u_2)
+\phi(a,u_1u_2,b_2)$$
$$~~~+\phi(u_1,a_1b_1,u_2)-\phi(u_1,a_1u_2,b_2)+
\phi(u_1,u_2a_2,b_2)-\phi(a,b,u_1u_2)
+\phi(a,u_1,b_1u_2)$$
$$~~~-\phi(a,u_1,u_2b_2)-\phi(u_1,a_1,b_1u_2)+\phi(u_1,a_1,u_2b_2)
-\phi(u_1,u_2,a_2b_2)
+\phi(a,b,u_1)$$
$$~~~-\phi(a,u_1,b_1)+\phi(a,u_1,u_2)+\phi(u_1,a_1,b_1)
-\phi(u_1,a_1,u_2)+\phi(u_1,u_2,a_2)$$}

We now add these two expressions.
Using the identities $au_1=u_1a_1$, $bu_1=u_1b_1$, $a_1u_2=u_2a_2$,
and $b_1u_2=u_2b_2$, cancelling and collecting terms, yields the expression

\small{$$[\mu\delta + \delta\mu](\phi)(a,b,u_1,u_2)
=\phi(a,u_1,u_2)-\phi(u_1,a_1,u_2)+\phi(u_1,u_2,a_2)
+\phi(b,u_1,u_2)
-\phi(u_1,b_1,u_2)$$
$$\hskip .5in ~+\phi(u_1,u_2,b_2)
-\phi(ab,u_1,u_2)+\phi(u_1,a_1b_1,u_2)-\phi(u_1,u_2,a_2b_2)$$}
This expression is exactly $e_1^*\theta+e_2^*\theta -e_{12}^*\theta$
applied to $\phi$, hence the proof is complete.
\end{proof}

The inverse transgression formula implies that
a 2-gerbe $\phi$ on an orbifold groupoid $\G$ induces a
1-gerbe $\theta(\phi)$ on the
associated inertia groupoid
$\wedge\G$. Furthermore, two equivalent (isomorphic) 2-gerbes induce
equivalent (isomorphic) 1-gerbes on
the inertia groupoid.

    Recall that there is an embedding
    $e: \G\rightarrow \wedge \G$ by $e(x)=1_x$ where $1_x$ is the
identity arrow. The image  $e(\G)$ is often
    referred as non-twisted sector and other components of $\wedge \G$ are called {\em twisted sectors}.

    \begin{cor}
If $\phi$ is a cocycle, then $e^*\theta(\phi)$ is a coboundary.
    \end{cor}

    \begin{proof}
         $e^*\theta(\phi)(u,v)=\theta(1,u,v).$
         Using the embedding $\lambda: \G\rightarrow \G^2$ given by
         $x\rightarrow (i_x, i_x)$, we can pull back
         $e^*_1\theta(\phi)+e^*_2\theta(\phi)-e^*_{12}\theta(\phi)$
         in cohomology.
         Note that
$$\lambda^*e^*_1\theta(\phi)=\lambda^*e^*_2\theta(\phi)=\lambda^*e^*_{12}\theta(\phi)
         =e^*\theta(\phi).$$
         Therefore, $e^*\theta(\phi)=\delta e^*\mu(\phi)$ is a coboundary.
This implies that restricted to the untwisted sector, our cocycle
$\theta(\phi)$ gives rise to a trivial cohomology class.
\end{proof}

%
%
%

\section{The Inverse Transgression in the Case of a Finite Group}

In the case when the original orbifold is $[*/G]$ where $G$ is a
finite group, the inverse transgression has a classical interpretation
in terms of shuffle products. Recall that $\wedge [*/G]$ can be thought
of in terms of $G$ with the conjugation action; this breaks up
into a disjoint union of orbits of the form $G/Z_G(g)$, indexed by
conjugacy classes. Each of these is in turn equivalent to
$[*/Z_G(g)]$; so we have a Morita equivalence
$\wedge [*/G] \cong \sqcup_{(g)} [*/Z_G(g)]$. Hence we can restrict
our attention to these components; in particular we would like
to describe each
$\theta_g: C^k(G, U(1))\to C^{k-1}(Z_G(g), U(1))$.
Now for a finite group $G$, the cochain complex $C^*(G,U(1))$ is in fact
equal to $Hom_G(B_*(G), U(1))$, where $B_*(G)$ is the bar resolution
for $G$ (see \cite{Brown}, page 19).

There is natural homomorphism $\rho_g: Z_G(g)\times\mathbb Z\to G$
given by $\rho_g(x, t^i) = xg^i$, where $t$ is a generator for $\mathbb Z$;
the
fact that $Z_G(g)$ centralizes $g$ is crucial here. This homomorphism
induces a map in integral homology

$$H_*(Z_G(g),\mathbb Z)\otimes H_*(\mathbb Z,\mathbb Z)\cong
H_*(Z_G(g)\times \mathbb Z, \mathbb Z)
\to H_*(G, \mathbb Z).$$
Classically it is known that multiplication
is induced by the shuffle product on the chain groups
(see \cite{Brown}, page 117--118);
i.e. there is a chain map
$B_*(Z_G(g))\otimes B_*(\mathbb Z)\to B_*(G)$ which will
induce $\rho_{g*}$ in homology.
Let $t$ denote a generator of the cyclic group
$\mathbb Z$. The shuffle product we are
interested in is
$B_k(Z_G(g))\otimes B_1(\mathbb Z)\to B_{k+1}(G)$
given by
$$[g_1|g_2|\dots | g_k]\star [t^i] =
\sum_{\sigma} \sigma [g_1|g_2|\dots |g_k|g_{k+1}]$$
where $g_{k+1}=g^i$,
$\sigma$ ranges over all $(k,1)$--shuffles and
$$\sigma [g_1|g_2|\dots |g_{k+1}] = (-1)^{sign(\sigma)}
[g_{\sigma(1)}|g_{\sigma(2)}|\dots |g_{\sigma(k+1)}].$$
A $(k,1)$--shuffle is an element $\sigma\in S_{k+1}$ such that
$\sigma (i) < \sigma (j)$ for $1\le i<j\le k$. These are precisely
the cycles:
$$1, (k~ k+1), (k-1~k~k+1), (k-2~k-1~k~k+1),\dots,
(1~ 2~ 3\dots k~ k+1).$$ Note that there are $k+1$ of them. This can
be dualized, using $U(1)$ coefficients, but for cohomology purposes
it's easier to use integral coefficients. Given a cocycle $\phi\in
C^{k+1}(G, \mathbb Z)$, we see that $\theta_g (\phi)\in
C^k(Z_G(g),\mathbb Z)$ can be defined as
$$\theta_g (\phi)([g_1|g_2|\dots |g_k]) = \phi([g_1|g_2|\dots | g_k]\star [g])$$
where $g_1,g_2,\dots , g_k\in Z_G(g)$.

As a consequence of this we see that
$\theta_g^*:H^{k+1}(G, \mathbb Z)\to H^k(Z_G(g), \mathbb Z)$
is induced by the multiplication map
$$\rho_g^*:H^{k+1}(G,\mathbb Z)\to H^k(Z_G(g), \mathbb Z)\otimes H^1(\mathbb Z, \mathbb Z).$$
To be precise, if $\nu$ is the natural generator for $H^1(\mathbb Z, \mathbb Z)$,
then
$$\rho_g^*(u) = res^G_{Z_G(g)}(u)\otimes 1 + \theta_g^*(u)\otimes\nu.$$
This discussion clarifies the geometric arguments in \cite{DW}, and will
also allow us to do some computations in cohomology.

\begin{example}\label{elementary abelian}
Finite group cohomology is difficult to compute, especially over
the integers. The simple examples such as cyclic and quaternion
groups are not so interesting in this context, as their odd dimensional
cohomology (with trivial $\mathbb Z$ coefficients) is zero.
The first interesting example is $G=(\mathbb Z/2)^2$. In this case
$H^*(G,\mathbb F_2)$ is a polynomial algebra on two degree one
generators $x, y$. In degree four there is a natural basis given by
$x^4, y^4, x^3y, x^2y^2, xy^3$. For an elementary abelian
$2$--group, the mod 2 reduction
map for $k>0$ is a monomorphism $H^k(G,\mathbb Z)\to H^k(G,\mathbb F_2)$,
and so we can understand it as the kernel of the Steenrod operation
$Sq^1: H^k(G,\mathbb F_2)\to H^{k+1}(G,\mathbb F_2)$. Hence we see
that $H^4(G,\mathbb Z)$ can be identified with the subspace generated
by $x^4$, $y^4$ and $x^2y^2$. These are all squares, hence when we
apply $\theta_g^*: H^4(G,\mathbb Z)\to H^3(G,\mathbb Z)$ for
any $g\in G$, the result
will always be zero.

Next we consider $G=(\mathbb Z/2)^3$; in this case $H^*(G,\mathbb F_2)$
is a polynomial algebra on three degree one generators $x,y,z$. In this
case we have an element $\alpha = Sq^1(xyz) = x^2yz + xy^2z + xyz^2$ which
represents
a non-square element in $H^4(G,\mathbb Z)$. By analyzing the multiplication
map in cohomology we obtain the following.

\begin{lem}
Let $g=x^ay^bz^c$ be an element in $G=(\mathbb Z/2)^3$, where we are
writing it in terms of the standard basis (identified with its dual
by abuse of notation). Then
$$\theta_g^*(\alpha) = a(y^2z+z^2y) + b(x^2z+xz^2) + c(x^2y+xy^2)$$
and so it is non--zero on every component except the one corresponding
to the trivial element in $G$.
\end{lem}

Now for an abelian group, the multiplicative formula implies
that for all $g,h\in G$,
$\theta_g^* + \theta_h^* = \theta_{gh}^*$ in cohomology,
or up to coboundaries. In particular
this implies that the correspondence
$g\mapsto \theta_g(\alpha)$ defines a homomorphism $G\to H^3(G,\mathbb Z)$
of elementary abelian groups, in this case an isomorphism.
\end{example}

\section{The Twisted Pontryagin Product for Finite Groups}

Let $G$ denote a finite group, and consider the orbifold defined
by its action on a point. Then the inertia groupoid
$\wedge \G$ can be identified with the groupoid determined by the
conjugation action of $G$ on itself.
In this case the untwisted orbifold K--theory is simply
$K_G(G)$, which is additively isomorphic to
$\bigoplus_{(g)}R(Z_G(g))$, where as before $Z_G(g))$ denotes
the centralizer of $g$ in $G$, and the sum is taken over
conjugacy classes. This group can be endowed with
a certain product, known as the \textsl{Pontryagin product},
defined as follows.
An equivariant vector bundle over $G$ (with the conjugation action)
can be thought of as a collection of finite dimensional vector spaces
$V_g$ with a $G$-module structure on $\oplus_{g\in G} V_g$ such that
$gV_h=V_{ghg^{-1}}$. The product of two of these bundles
is now defined as:
$$(V\star W)_g = \bigoplus_{\{g_1,g_2\in G, ~~g_1g_2=g\}}
V_{g_1}\otimes W_{g_2}$$
This formula has been referred to as
the holomorphic orbifold model in the physics literature \cite{DVVV}.

This product admits an alternate description, which will admit
a geometric generalization.
In this case we can identify $\G^2$ with the orbifold defined by
considering $G\times G$ with the conjugation action on both
coordinates. Our maps $e_1$, $e_2$ and $e_{12}$ correspond to
$(g,h)\mapsto g$, $(g,h)\mapsto h$, $(g,h)\mapsto gh$
respectively, which are $G$--equivariant with respect to
the conjugation action.
Then, if $\alpha$, $\beta$ are elements in $K_G(G)$,
the Pontryagin product can also be defined as
$$\alpha\star\beta =
e_{12*}(e_1^*(\alpha)\cdot e_2^*(\beta)).$$

We propose to extend this definition to twisted K--theory,
with certain conditions on the twisting cocycle.
Note that given a $2$--cocycle
$\tau=\theta (\phi)$
in the image
of the inverse transgression, then by our multiplicative
formula we have
$e_1^*\tau + e_2^*\tau = e_{12}^*\tau + \delta\mu (\phi)$.

\begin{dfn}
Let $\tau$ be a $U(1)$ valued $2$--cocycle for the orbifold
defined by the conjugation action of a finite group $G$
on itself which is in the image of the inverse transgression.
The Pontryagin product on $^\tau K_G(G)$ is defined by the following
formula: if $\alpha, \beta \in ~^\tau K_G(G)$, then

$$\alpha\star\beta = e_{12*}(e_1^*\alpha\cdot e_2^*\beta)$$

\end{dfn}
Note that
if $\tau=\theta (\phi)$
then by our multiplicative
formula we have
$$e_1^*\tau + e_2^*\tau = e_{12}^*\tau + \delta\mu (\phi)$$
and so
the product $e_1^*(\alpha)\cdot e_2^*(\beta)$ lies in
$$^{e_1^*\tau + e_2^*\tau}K_G(G)=
~^{e_{12}^*\tau +
\delta\mu (\phi)}K_G(G)\cong ~^{e_{12}^*\tau}K_G(G).$$
Now applying $e_{12*}$, this is mapped to $^\tau K_G(G)$;
and so we have a product on our twisted K--theory;
it is elementary to verify that this defines an associative
product.

Our approach will define a twisted Pontryagin
product for any cocycle in the image of the inverse
transgression.
This cocycle could very well be a coboundary; but that
does not necessarily imply that the corresponding
product on $K_G(G)$ is the standard Pontryagin
product. It is also clear that we may choose twistings
which give rise to a twisted K--theory without any
product.\footnote{For much more on Pontryagin products please
consult \cite{FHT} and its sequels.}

\begin{example}
If $G$ is an abelian group, then what we are doing
is using the identification
$\theta (\phi)_g + \theta (\phi)_h = \theta (\phi)_{gh}$
to define a product on the abelian group
$$\mathcal{X}(G)= \sum_{g\in G} ~ ^{\theta (\phi)_g} R (G)$$
via the pairing
$$^{\theta(\phi)_g} R(G)\otimes~^{\theta(\phi)_h}R(G)
\to~^{\theta (\phi)_{gh}}R(G).$$

In the case described in \ref{elementary abelian}
for $G=(\mathbb Z/2)^3$ and a cocycle
$\phi$
representing the cohomology
class $xy^2z + xyz^2 + x^2yz$, $\theta (\alpha)$
establishes a group homomorphism
$G\to H^3(G,\mathbb Z)$, with image the
subgroup generated by $xy^2+x^2y$, $xz^2+x^2z$ and
$yz^2 + y^2z$.
In this case we see that for $g\ne 1$,
$^{\theta (\phi)_g}R(G)$ has rank equal to two, and so
$\mathcal{X}(G)$ is of rank equal to twenty-two, with
the twisted Pontryagin product described above.
This can be made explicit.
\end{example}

The case of the Pontryagin product should be considered
as motivation for the case of orbifold groupoids. As long
as we twist with a cocycle in the image of the inverse
transgression, the levels will match up as required. Hence
the main difficulty is geometric--as we shall see in the next
section, there is an obstruction bundle which plays an important
role.

\section{Twisted K--theory of Orbifolds}

During the course of our investigation of possible stringy products on the
twisted K-theory of
orbifolds, we came to realize that one needs to use the
very same
information required
to construct the Chen-Ruan cohomology of orbifolds (\cite{CR}).
We first briefly recall the situation for orbifold
cohomology, and then proceed to develop the tools necessary to
deal with twisted K--theory. For a very interesting but
different approach we refer the reader to \cite{JKK}.

First we recall that $H^*_{CR}(\G,\mathbb C)$ is additively
the same as $H^*(\wedge\G, \mathbb C)$; what is interesting
is the ring structure.
Recall that there are
three evaluation maps $e_1, e_2, e_{12}: \G^2\to
\wedge\G$. A naive definition of the stringy product
for orbifold cohomology
(which would generalize the
Pontryagin product)
would be $\alpha\star \beta= (e_{12})_*(e^*_1\alpha\cup e^*_2 \beta)$.
However,
one soon discovers that this product is not associative due to the
fact that $e_1, e_{12}$ are not transverse
in general. In fact, the correction term is precisely the obstruction bundle
to the transversality of
$e_1$ and $e_{12}$. A natural idea was to modify the definition of the product
via
$$\alpha\star \beta=(e_{12})_*(e^*_1\alpha\cup e^*_2 \beta\cup e(E_{\G^2}))$$
where we need to construct a bundle $E_{\G^2}$ in a fashion that is
consistent with the obstruction to transversality of $e_1, e_{12}$,
and $e(E_{\G^2})$ denotes its Euler class. A key observation is that
the \textsl{obstruction bundle} in the construction of the Chen-Ruan
product provides such a choice. We will adapt this same idea to
$K$--theory.


Throughout this section, we assume that $\G$ is a compact, almost
complex orbifold groupoid. Then $\G^k$ naturally inherits an almost
complex structure such that the evaluation map $e_{i_1, \cdots,
i_l}: \G^k\rightarrow \G^l$ is an almost complex quasi-embedding.
$\G^2$ can be identified with the space of constant morphisms
     from an orbifold sphere with
    three marked point to $\G$; we now make this identification
    precise. Consider an orbifold Riemann sphere with three
    orbifold points
   $(\mathbb S^2, (x_1, x_2, x_3), (m_1, m_2, m_3))$; in this context we simply denote it by
   $\mathbb S^2$.
   Suppose that $f$ is a constant morphism from $\mathbb S^2$ to $\G$. Here, the term constant means that the induced
   map $|f|: \mathbb S^2\rightarrow |\G|$ is constant. Let $y=\rm{im}(|f|)$ and $U_y/G_y$ be an orbifold chart at $y$.
   By the results in \cite{ALR}, $f$ is classified by the conjugacy class of a homomorphism
   $\pi_f: \pi^{orb}_1(\mathbb S^2)\rightarrow G_y$.
    Recall that
   $$\pi^{orb}_1(\mathbb S^2)=\{\lambda_1, \lambda_2, \lambda_3;
   \lambda^{k_i}_i=1, \lambda_1\lambda_2\lambda_3=1\},$$
   where $\lambda_i$ is represented by a loop around the marked point
   $x_i$. $\pi_f$ is uniquely determined by a pair of elements $(g_1, g_2)$ with $g_i\in G_y$
   where $g_i=\pi_f(\lambda_i)$; on the other hand, $(g_1,g_2)\in \G^2_0$. It is clear that the same method
   can be used to identify the moduli space of constant morphisms from an orbifold sphere with $k+1$ marked points
   to $\G$ with the groupoid of k--multisectors
   $\G^k$.

   For any $f\in \G^2$ viewed as a constant morphism, we can form an elliptic complex
   $$\bar{\partial}_f: \Omega^0(f^*T\G)\rightarrow \Omega^{0,1}(f^*T\G),$$
   where $f^*T\G$ is a complex vector bundle by our assumption.
   This defines an orbibundle $E$ with $E_{\G^2}|_{f}=\rm{Coker}~\bar{\partial}_f$.
   We now examine $E_{\G^2}$ in more detail.
      Let $g_1,g_2\in \G^2_0$; by the definition, $g_1,g_2\in G_x$ for $x=s(g_i)=t(g_i)$.
      Let $N$ be the subgroup of $G_x$ generated by $g_1, g_2$.
   By Lemma 4.5 in \cite{ALR}, $N$ depends only on the component of $\G^2$.  Let
   $e: \G^2\rightarrow \G$ be an evaluation map. Clearly
   $N$ acts on $e^*T\G$  while fixing $T\G^2$.

       There is an obvious surjective homomorphism
   $\pi: \pi^{orb}_1(\mathbb S^2)\rightarrow N$
   and $\rm{Ker}~\pi$ is therefore a subgroup of finite index. Suppose that
   $\tilde{\Sigma}$ is the orbifold universal cover of $\mathbb S^2$.
   By \cite{ALR} (see Chapter II),  $\tilde{\Sigma}$
   is smooth. Let $\Sigma=\tilde{\Sigma}/\rm{Ker}~\pi$; then $\Sigma$ is
   compact and we have a quotient map
   $\Sigma\rightarrow \mathbb S^2=\Sigma/N$. Since $N$ contains the relation
   $g^{m_i}_i=1$, $\Sigma$ is smooth.

    It is clear that $f$ lifts to an ordinary constant map
    $\tilde{f}: \Sigma\rightarrow U_y$;
hence $\tilde{f}^*T\G=T_y \G$ is a trivial bundle over $\Sigma.$
Then we can lift the elliptic complex to $\Sigma$
$$\bar{\partial}_{\Sigma}: \Omega^0(\tilde{f}^*T\G)\rightarrow \Omega^{0,1}(\tilde{f}^*T\G).$$
The original elliptic complex is just the $N$-invariant part of
the current one. However,  $\rm{Ker}(\bar{\partial}_{\Sigma})=T_y \G$
and $\rm{Coker}~(\bar{\partial}_{\Sigma})=H^{0,1}(\Sigma)\otimes T_y
\G$. Now we vary $y$ in a component $\G^2_{(\g)}$ to obtain
$e^*_{(\g)}T\G$ for the evaluation map
     $e_{(\g)}: \G^2_{(\g)}\rightarrow \G$ and $H^{0,1}(\Sigma)\otimes e^*_{(\g)}T\G$. $N$ acts on both.
    It is clear that $(e^*_{(\g)}T\G)^N=T \G_{(\g)}$ as we claim.   The obstruction bundle $E_{(\g)}$ we want is
   the invariant part of $H^{0,1}(\Sigma)\otimes_{\mathbb C} e^*_{(\g)}T\G$, i.e.,
   $E_{(\g)}=(H^{0,1}(\Sigma)\otimes_{\mathbb C} e^*_{(\g)}T\G)^N$. We remark that $E_{(\g)}$
   can have different dimensions
   at the different components of $\G^2$.
   We can obviously use the same method to construct a bundle $E_{\G^k}$
   over $\G^k$ whose fiber is the
   cokernel of $\bar{\partial}_f$ for $f\in \G^k$.

   Consider the evaluation maps $e_1,e_2,e_{12}: \G^2\rightarrow \wedge\G$.
    Suppose that the local chart of $\G$ is $U/G$. Then, the
   local chart of $\wedge\G$ is $(\sqcup_{g\in G}U^g)/G$.
    The local chart of $\G^2$ is $(\sqcup_{g_1,g_2\in G} U^{g_1}\cap U^{g_2})/G$. 
    The $e_1,e_2,e_{12}$ are quasi-embeddings
   and hence $\G^2$ can be thought of as a quasi--suborbifold of $\wedge\G$ via three
    different quasi-embeddings, denoted by $e_{1*}\G^2, e_{2*}\G^2, e_{12*}\G^2$.
   It is clear that
   $\G^3=\G^2 \times_{e_{12},e_1} \G^2$ and so $\G^3$ can be viewed as the intersection
   of the quasi--suborbifolds $e_{12*}\G^2$ and $e_{1*}\G^2$:

$$
\begin {array}{rcl}
\G^3& \stackrel {\pi_2}{\to}&\G^2\\
 \pi_1\downarrow&& \downarrow e_1 \\
 \G^2& \stackrel {e_{12}}{\to}&\wedge\G
\end{array}
$$

   The problem is that $e_{12*}\G^2, e_{1*}\G^2$
   are not intersecting
   transversely in general. Let
   $$\nu=(e_{12}\pi_1)^*T\wedge\G/\pi_1^*T\G^2+\pi^*_2T\G^2$$
   where $\pi_i: \G^3\rightarrow \G^2$ are the natural projection maps. This is the so-called
   excess bundle of the intersection. A crucial ingredient in the proof of
associativity for the Chen-Ruan
   product is the bundle formula:
   \begin{Th}
   $$E_{\G^3}\cong \pi^*_1 E_{\G^2}\oplus \pi^*_2 E_{\G^2}\oplus \nu.$$
   \end{Th}
\noindent The proof is given in \cite{CR} by gluing arguments.
   We shall call this the \textsl{obstruction bundle equation}.

   These bundles will play a key role in our definition of the
   product. It ties in to the transversality question mentioned
   before via the following lemma, which is an orbifold analogue
   of the clean intersection formula first described by Quillen
   \cite{Q}:

   \begin{lem}
   Suppose that $i_1:\H_1\to \G$ and $i_2:\H_2\to\G$ are
   quasi--suborbifolds of $\G$ forming a clean
   intersection $\H_3$. Then, if $u\in~
   ^\alpha K(\H_1)$ we have
   $$i_2^*{i_1}_*u={\pi_2}_*[\pi_1^*u\cdot e_K(\nu)]$$
   where $\pi_1:\H_3\to \H_1$ and $\pi_2:\H_3\to \H_2$ are the natural
   projections, $\nu$ is the excess bundle of the intersection,
   and $e_K(\nu)\in K(\H_3)$ denotes
   its Euler class.
   \end{lem}

   Note that this lemma will allow us to connect facts about the
   geometry of the quasi--suborbifolds $\G^2$ and $\G^3$ in $\wedge\G$
   with products, and the obstruction bundle $E$ plays a key role
   here. In fact the only information we need about $E$ is the
   obstruction bundle equation mentioned above.

   \begin{dfn}
   Suppose that $\G$ is an almost complex orbifold.
Let $\phi$ denote a $U(1)$--valued $3$--cocycle for $\G$, and $\theta(\phi)$ its
inverse transgression, i.e. a $U(1)$--valued $2$--cocycle for $\wedge\G$.
For $\alpha, \beta\in ~^{\theta (\phi)} K(\wedge\G)$, we define
   $$\alpha\star \beta  =
   e_{12*}(e^*_1\alpha\cdot e^*_2\beta\cdot e_K(E_{\G^2})).$$
   \end{dfn}

   Our goal is to show that this defines an
   associative product. However, we want to do this using very
   general properties of our construction, so that it will be
   a natural extension of the Chen--Ruan product. We will
   abbreviate $E_2 = E_{\G^2}$, $E_3= E_{\G^3}$.\footnote{We have chosen to
   work only with even dimensional $K$--theory, but this can be
   readily extended to odd dimensions; the signs work out
   appropriately after a tedious computation.}


As with the twisted Pontryagin
product, we must explain in what sense this defines a product.
Given that $e_1,e_2, e_{12}: \G^2\to\wedge G$,
then $e_i^*(\alpha)\in~ ^{e_i*(\theta (\phi))}K(\G^2)$
for $i=1,2$. Hence

$$e_1^*(\alpha)\cdot
                   e_2^*(\beta)\cdot e_K(E_{\G^2})
\in ~^{e_1^*\theta (\phi) + e_2^*\theta (\phi)}K(\G^2).$$
Now here we must be careful. We use the fact that the twisting
cocycle is in fact equal to $e_{12}^*\theta(\phi) + \delta\mu(\phi)$.
We have a \textsl{canonical} isomorphism

$$^{e_{12}^*\theta(\phi) + \delta\mu(\phi)} K(\G^2)
\cong ~^{e_{12}^*\theta(\phi)} K(\G^2).$$
Next we apply the push forward $e_{12*}$ to our expression,
which now lands in
$^{\theta(\phi)}K(\wedge \G)$. Note that the product is clearly
commutative; it is associativity that requires a proof.

   \begin{Th}
   This product is associative:
   $$(\alpha\star \beta)\star\gamma = \alpha\star(\beta\star \gamma).$$
   \end{Th}

   \begin{proof}

We pull back the quasi-embeddings $e_1$, $e_2$ and $e_{12}$ to
$\G^3$ as follows: let $\tilde{e_1}=e_1\pi_1$;
$\tilde{e_2}=e_2\pi_1$; $\tilde{e}_3=e_2\pi_2$ and
$\tilde{e}_{123}=e_{12}\pi_2$. We will also make use of an operator
defined for triples, namely let $I_3: \G^3\to\G^3$ be defined by
$$(g_1,g_2,g_3)\mapsto (g_2, g_3, g_3^{-1}g_2^{-1}g_1g_2g_3).$$
Note that $\tilde{e_i}I_3$ is equal to $\tilde{e}_{i+1}$ up to
conjugation (modulo three);
hence they induce the same map in K-theory. Note that
$e_K(E_3)$ is invariant under $I_3^*$ by construction.

We are now ready to start computing:

$$(\alpha\star\beta)\star\gamma =    e_{12*}(e_1^*(\alpha\star\beta)\cdot e_2^*\gamma\cdot
    e_K(E_2))
    =
    e_{12*}(e_1^*(e_{12*}(e_1^*\alpha\cdot e_2^*\beta\cdot
    e_K(E_2)))\cdot e_2^*\gamma\cdot e_K(E_2))$$
    We can apply the clean intersection formula to this
    expression, whence we obtain

    $$(\alpha\star\beta)\star\gamma  =  e_{12*}(\pi_{2*}(\pi_1^*(e_1^*\alpha\cdot
e_2^*\beta\cdot
    e_K(E_2))\cdot e(\nu))
    \cdot e_2^*\gamma\cdot e_K(E_2))$$
    $$=e_{12*}(\pi_{2*}(\tilde{e_1}^*\alpha\cdot\tilde{e}_2^*\beta
    \cdot\pi_1^*e_K(E_2)\cdot e_K(\nu))\cdot e_2^*\gamma\cdot e_K(E_2))$$
    $$=e_{12*}\pi_{2*}(\tilde{e}_1^*\alpha\cdot\tilde{e}_2^*\beta
    \cdot\pi_1^*e_K(E_2)\cdot
    e_K(\nu)\cdot\pi_2^*e_2^*\gamma\cdot\pi_2^*e_K(E_2))$$
    Note here that we are using the following property of the
    pushforward: $\pi_{2*}(x\cdot \pi_2^*y)=\pi_{2*}(x)\cdot y$.
    Thus we have
    $$(\alpha\star\beta)\star\gamma  =  \tilde{e}_{123*}(\tilde{e}_1^*\alpha
    \cdot\tilde{e}_2^*\beta\cdot\tilde{e}_3^*\gamma
    \cdot\pi_1^*e_K(E_2)\cdot\pi_2^*e(E_2)\cdot e_K(\nu))$$
    $$=\tilde{e}_{123*}(\tilde{e}_1^*\alpha\cdot
    \tilde{e}_2^*\beta\cdot \tilde{e}_3^*\gamma\cdot e_K(E_3))$$

    Now we consider
    $$\alpha\star (\beta\star\gamma) =
    (\beta\star\gamma)\star\alpha
    =
    \tilde{e}_{123*}(\tilde{e}_1^*\beta\cdot\tilde{e}_2^*\gamma
    \cdot\tilde{e}_3^*\alpha\cdot e_K(E_3))$$
    $$=\tilde{e}_{123*}I_{3*}I_3^*
    (\tilde{e}_1^*\beta\cdot\tilde{e}_2^*\gamma\cdot
    \tilde{e}_3^*\alpha\cdot e_K(E_3))$$
    $$=\tilde{e}_{123*}(I_3^*\tilde{e}_1^*\beta
    \cdot I_3^*\tilde{e}_2^*\gamma\cdot I_3^*\tilde{e}_3^*\alpha
    \cdot I_3^*e_K(E_3))$$
    $$=\tilde{e}_{123*}(\tilde{e}_2^*\beta\cdot \tilde{e}_3^*\gamma
    \cdot\tilde{e}_1^*\alpha\cdot e_K(E_3))$$
    $$=
    \tilde{e}_{123*}(\tilde{e}_1^*\alpha\cdot\tilde{e}_2^*\beta
    \cdot\tilde{e}_3^*\gamma\cdot e_K(E_3))$$
$$=(\alpha\star\beta)\star\gamma$$
Hence our proof of associativity is
    complete.
\end{proof}

\end{document}